\documentclass[12pt]{article}

\usepackage{amsmath}
\usepackage{amssymb}

\usepackage[T1]{fontenc} 
\usepackage{sets pace}
\usepackage{comment}
\usepackage{colortbl}
\definecolor{text1}{cmyk}{1,.65,0,0} 
\definecolor{text2}{rgb}{1,0,0} 
\definecolor{text3}{cmyk}{0,0,0,1} 
\definecolor{text4}{cmyk}{0,0,0,0.5} 
\definecolor{text5}{cmyk}{1.0,0.0,1.0,0} 

\usepackage{amsfonts}
\usepackage{amsthm}
\usepackage{graphicx}
\usepackage{caption}
\usepackage{subfigure}
\usepackage{latexsym}
\usepackage[most]{tcolorbox}
\makeatletter
\renewcommand{\@seccntformat}[1]
{\csname the#1\endcsname.\enspace}
\makeatother
\newtheorem{definition}{Definition}
\newtheorem{theorem}{Theorem}
\newtheorem{lemma}{Lemma}
\newtheorem{remark}{Remark}
\newtheorem{corollary}{Corollary}

\usepackage{cancel}
\usepackage[normalem]{ulem}

\usepackage{natbib}

\begin{document}
\large	
\begin{center}
	\textbf{Bayesian prediction regions and density estimation with type-2 censored data}  \footnote{\today. }\\
\end{center}	
\normalsize
\begin{center}
	{\sc Akbar Asgharzadeh$^{a}$, \'Eric Marchand$^{b}$ \& Ali Saadati Nik$^{a}$ } \\
	{\it a  Department of Statistics, University of Mazandaran, P.O. Box 47146-1407, Balbosar, IRAN \\b  D\'epartement de math\'ematiques, Universit\'e de
		Sherbrooke, Sherbrooke Qc,
		CANADA, J1K 2R1 \\
		\quad (e-mails:  a.asgharzadeh@umz.ac.ir;\, eric.marchand@usherbrooke.ca; \,  a.saadatinik@umz.ac.ir) } \\
\end{center}
\begin{center}
	{\sc Summary}
\end{center}
\small  
For exponentially distributed lifetimes, we consider the prediction of future order statistics based on having observed the first $m$ order statistics.  We focus on the previously less explored aspects of predicting: (i) an arbitrary pair of future order statistics such as the next and last ones, as well as (ii) the next $N$ future order statistics.   We provide explicit and exact Bayesian credible regions associated with Gamma priors, and constructed by identifying a region with a given credibility $1-\lambda$ under 
the Bayesian predictive density.  For (ii), the HPD region is obtained, while a two-step algorithm is given for (i).   The predictive distributions are represented as mixtures of bivariate Pareto distributions, as well as multivariate Pareto distributions.    For the non-informative prior density choice, we demonstrate that a resulting Bayesian credible region has matching frequentist coverage probability, and that the resulting predictive density possesses the optimality properties of best invariance and minimaxity.

\noindent  
\vspace{1.5mm}

\vspace{1.5mm}


\section{Introduction}

Predictive analysis  based on censored data in life testing experiments is fundamental and leads to interesting challenges.   In this work, we are concerned with prediction regions for future order statistics based on the first $m$ order statistics generated by exponentially distributed data.   There has been some previous work on such problems, but we focus here on the less explored: (a) multivariate aspects, and (b) use of Bayesian predictive densities to generate prediction regions and their related theoretical properties.   

We consider an i.i.d. sample of size $n$ from an exponential distribution with density $\theta e^{-\theta t} \mathbb{I}_{(0,\infty)}(t)$, but we are only able to observe the first $m$ order statistics $X_{1:n}, \ldots, X_{m:n}$, commonly referred to as type-II censoring.   Our objectives include:  
\begin{enumerate}
\item[\bf (I)]   the joint prediction of two future order statistics $X_{r:n}$ and $X_{s:n}$, with $m<r<s \leq n$, as well as

\item[\bf (II)]  the joint prediction of the next $N$ order statistics with $2 \leq N \leq n-m$.  

\end{enumerate}
Scenario {\bf (II)} seems a natural one to consider and includes the case of all future order statistics with $N=n-m$.   Scenario {\bf (I)} in more relevant in situations where the focus is for instance on the next and last order statistics (i.e., $r=m+1, s=n$), or the last two order statistics (i.e., $r=n-1, s=n$).   Obviously, both scenarios overlap for $N=2$.    Surprisingly perhaps, for Gamma priors which  we consider, the specification of Bayesian predictive densities and regions leads to more complex representations in the bivariate case.  The bivariate case {\bf (I)} was considered recently by \cite{Bagheri-et-al} and we refer to this work for further motivation and historical aspects of the problem. Their prediction regions are non-Bayesian however and based on a pair pivotal quantities which will not arise as a Bayesian solution. 

It seems natural to generate prediction regions for future order statistics via a Bayesian predictive density, but such problems seem to have been relatively unexplored.  An exception is given by \cite{Dunsmore-1974}, but his work concerns a single order statistic. We proceed in doing so for Gamma prior densities as well as the non-informative density $\pi_0(\theta)=\frac{1}{\theta} \, \mathbb{I}_{(0,\infty)}(\theta)$.    The obtained expressions are quite tractable and, interestingly,  bring into play finite mixtures of bivariate Pareto distributions, and multivariate Pareto distributions.   The obtained mixtures are ``non-probabilistic'' in the sense that the coefficients $a_i$ in $\sum_i a_i f_i(t)$; the $f_i$'s being densities; take on both positive and negative values.   This aspect is less appreciated than ``probabilistic mixtures'', but it does not hinder the usefulness of the representation for computational purposes of moments and cumulative distributions functions.

For a given credibility $1-\lambda$, we fully describe the highest posterior density (HPD) prediction region for scenario {\bf (II)}, while we propose an algorithm to generate exact prediction regions for the bivariate scenario {\bf (I)}, based on a natural decomposition of the joint predictive density into its marginal and conditional parts.    The analysis is carried out and much facilitated by considering the prediction of spacings between future order statistics, which can be converted back to the prediction of the order statistics themselves.

Many of the recent studies on predictive density estimation focus on theoretical properties of the density itself in a decision-theoretic context.   For scenarios {\bf (I)} and {\bf (II)}, we report on the best invariant and minimax properties for Kullback-Leibler divergence loss of the Bayesian predictive density associated with the usual non-informative prior density $\pi(\theta)=\frac{1}{\theta}$, attributable to the work of \cite{Liang-Barron}.  Moreover, we show that the frequentist probabilty of coverage for a prediction region generated by such a Bayesian predictive density matches its exact credibility.   Therefore, such a Bayesian density not only possesses optimality properties in a decision-theoretic framework, but also provides a satisfactory frequentist option that compares favourably with previous solutions (e.g., Bagheri et al., 2022).

The organization of the manuscript is as follows.    The preliminary results of Chapter 2 cover model densities for scenarios {\bf (I)} and {\bf (II)}, some general aspects on Bayesian predictive densities, and multivariate Pareto distributions.  Bayesian predictive densities and regions are derived and illustrated in Section 3.  In Section 4, we demonstrate in a more general context including ours that a Bayesian credible region with credibility $1-\lambda$ associated with the non-informative prior density $\pi_0$ yields matching frequentist coverage probability $1-\lambda$ for all $\theta>0$.    Finally in Section 5, we concisely review optimality properties of the Bayesian predictive densities again associated with density $\pi_0$ and pertaining to the best invariant and minimaxity properties.

\section{Preliminary results}

\subsection{Model densities}

For the first $m$ order statistics $X_{1:n}, \ldots, X_{m:n}$ among $n$ generated from i.i.d. $Exp(\theta)$ data, it is well known (e.g.,  \cite{Lawless-1971}) $X=\sum_{i=1}^m X_{i:n} + (n-m) X_{m:n}$ is a sufficient statistic and Gamma distributed $\mathcal{G}(m,\theta)$.  Hereafter, we therefore consider such a summary and corresponding density, that is
\begin{equation}
\label{modelsufficient}
X \sim  \frac{\theta^m}{\Gamma(m)} \, x^{m-1} \, e^{-\theta x} \, \mathbb{I}_{(0,\infty)}(x).
\end{equation}

For the prediction of $X_{r:n}$ and $X_{s:n}$ based on $X$, it is convenient and equivalent to consider 
\begin{equation}
\label{predicted}
Y=(Y_1, Y_2)^{\top} \hbox{ with }  Y_1\,=\,  X_{r:n} - X_{m:n}, \, Y_2\,=\,   X_{s:n} - X_{r:n}.
\end{equation}
The equivalence stems from the correspondence between a prediction region $R$ for $Y$ and the inverse mapping $\{(y_1+x_{m:n}, y_1+y_2+x_{m:n}): (y_1, y_2) \in R \}$ as a prediction region for $(X_{r:n}, X_{s:n})$.

For the multivariate version where the objective is to predict jointly the $N$ next future order statistics, it is analogously useful to consider
\begin{equation}
\label{predictedmultivariate}
Z=(Z_1, \ldots, Z_N)\,, \hbox{ with } \, Z_i=X_{m+i:n}-X_{m+i-1:n} \hbox{ for } i \in \{1, \ldots, N\},   2 \leq N \leq n-m\,,
\end{equation}
as the objects of prediction.
 
As described with the next result, both transformations (\ref{predicted}) and (\ref{predictedmultivariate}) lead to convenient underlying model densities for $Y$ or $Z$.

\begin{lemma}
\label{lemmadensities}
Given set-up (\ref{modelsufficient}, \ref{predicted}, and \ref{predictedmultivariate}) and fixed $\theta$:
\begin{enumerate}
\item[\bf{ (a)}]   $Y$ and $X$ are independently distributed; 
\item[\bf{ (b)}]  $Z$ and $X$ are independently
  distributed;
\item[\bf{ (c)}]   $Y_1$ and $Y_2$ are independently distributed with joint density on $\mathbb{R}_+^2$ given by $q_{\theta}(y) \, = \, \theta^2 \,  q_1(\theta y_1, \theta y_2)$ and with 
\begin{equation}
\label{lemmadensityofY}
q_1(y) \, = \, \frac{(n-m)!}{(r-m-1)! \, (s-r-1)! \, (n-s)!} \, 
\big(1-e^{-y_1}\big)^{r-m-1} \, \big(1-e^{-y_2}\big)^{s-r-1} \, e^{-(n-r+1) y_1 - (n-s+1) y_2 }\,;
\end{equation}
\item[\bf{ (d)}]   $Z_1, \ldots, Z_N$ are independently distributed
   with  $Z_i \sim Exp\big((n-m-i+1) \, \theta\big)$.
\end{enumerate}

   \end{lemma}
{\bf Proof.}   
Part {\bf (b)} follows from a familiar renewal property of the exponential distribution, and furthermore implies {\bf (a)} since $Y$ is a function $h(Z)$ of $Z$.   For {\bf (c)}, since $(U=X_{r:n}-X_{m:n},V=X_{s:n}-X_{m:n})$ and $X_{m:n}$ are independently distributed, the distribution of $(U,V)$ matches that of its conditional distribution given $X_{m:n}$, and the latter joint distribution can by seen to be equivalent, with the above-mentioned renewal property, to that of the $(r-m)^{\hbox{th}}$ and $(s-m)^{\hbox{th}}$ order statistics from a sample of size $n-m$ from an $Exp(\theta)$ distribution.   Such a joint density is given by
\begin{equation}
f(u,v)  \, = \, \frac{(n-m)!}{(r-m-1)! \, (s-r-1)! \, (n-s)!} \, 
\big(1-e^{-\theta u}\big)^{r-m-1} \, \big(e^{-\theta u}-e^{-\theta v}\big)^{s-r-1} \, (e^{-\theta v})^{n-s} \, \theta^2 \, e^{-\theta (u+v)},\,   
\end{equation}
for $0 < u < v$.  The result follows by transforming $(U,V)$ to $Y$.
Finally for {\bf (d)}, since $Z$ and $X_{m:n}$ are independently distributed, the distribution of $Z$ matches that of its conditional distribution given $X_{m:n}$, and as above the latter joint distribution matches that of the distribution of the first $N$ order statistics spacings  from a sample of size $N$ from an $Exp(\theta)$ distribution.   It is well known that such order statistics spacings are independently distributed as $Exp\big((n-m-i+1)\theta\big)$ (e.g., Lehmann \& Casella, 1998, problem 6.18, page 71) which yields the result.  
\qed

\begin{remark}  The marginal densities of $Y_1$ and $Y_2$ in (\ref{lemmadensityofY}) can be extracted from part {\bf (c)} and are of the form $\mathcal{B}^{-1}(c_1,c_2) \, (1-e^{-u})^{c_1-1} \, e^{-c_2 u} \, \mathbb{I}_{(0,\infty)}(u)$, with  $\mathcal{B}(c_1,c_2)$ the Beta function.  An alternative and readily verified representation for the distribution of the above $Y_1$ and $Y_2$ is:
\begin{equation}
\nonumber
e^{- \theta Y_1} \sim \hbox{Beta}(n-r+1, r-m)\,,\, e^{-\theta Y_2} \sim \hbox{Beta}(n-s+1, s-r)\, \hbox{ independent}\,.  
\end{equation} 
\end{remark}
\subsection{Bayesian predictive densities}

 A general set-up for predictive density estimation relates to the following model density:
\begin{equation*}
\label{model2}  (X,Y)|\theta \sim p_{\theta}(x) \, q_{\theta}(y|x)\,, x \in \mathbb{R}^{d_1},y \in \mathbb{R}^{d_2}.
\end{equation*} 
We observe $X$ according to $p_{\theta}$ and we wish to estimate the density
 $q_{\theta}(\cdot|x)$.  Except for $\theta$, the densities are known.   The observed $X$ provides information about $\theta$ and determines the conditional density $q_{\theta}(\cdot|x)$ to estimate.   Much, but not all (e.g.,  Fourdrinier et al. 2019), previous work on properties of predictive densities focusses on models for which $X$ and $Y$ are conditionally independent, i.e.,  $q_{\theta}(\cdot|x) \equiv q_{\theta}(\cdot)$ for all $x$, but it is natural in general to estimate the density $q_{\theta}(\cdot|x)$ when there is dependence.   
 
Assume that we have a prior density $\pi$ for $\theta$, and a resulting posterior density $\pi(\cdot|x)$ taken to be absolutely continuous with respect to measure $\nu$.  A natural estimator of the density $q_{\theta}(\cdot|x)$ is the Bayes predictive density given by the conditional density $q(\cdot|x)$ of $Y|X=x$.  Integrating out $\theta$, we obtain the estimator
\begin{equation}
\label{predictivedensitygeneral}
\hat{q}_{\pi}(y|x)\,=\, \int_{\Theta} q_{\theta}(y|x) \, \pi(\theta|x) \, d\nu(\theta).
\end{equation}  
The above is a fully Bayesian procedure that can be used for obtaining predictive point estimates or prediction regions for $Y$.  The estimator or density $\hat{q}_{\pi}(\cdot;X)$ will naturally depend on $X$ through a sufficient statistics $T(X)$.

\subsection{Multivariate Pareto densities}
The predictive densities that we elicit below bring into play univariate, bivariate, and multivariate type II Pareto distributions, as well as non-probabilistic mixtures of such distributions.  Such distributions have been extensively studied (see for instance Section 52.4 of Kotz et al. 2000, or Arnold 2014, and the references therein).  
With convenient expressions for the moments and cumulative distribution functions (c.d.f.'s) relative to the $f_i$'s, the mixture representations clearly facilitate expressions for moments and c.d.f.'s for the full distribution.  

Multivariate Pareto distributions possess Pareto marginals with parameters $\ell,h>0$, densities and survival functions, which we denote and write as 
\begin{equation}
\label{densityunivariatePareto}
f_{\ell,h}(t) \, = \,  \frac{\ell\,h}{\;\;\;(1+ht)^{\ell+1}} \, \mathbb{I}_{(0,\infty)}(t)\,, \hbox{ and } \bar{F}_{l,h}(t) \, = \, (1+ht)^{-\ell}\,,
\end{equation} 
 for $t>0$, respectively.  Here is a multivariate Pareto definition; which includes the bivariate case; and some useful properties. 
\begin{definition}
\label{defmultivariatePareto}
A random vector $Z=(Z_1, \ldots, Z_N)^{\top}$ has a multivariate Pareto type II distribution, denoted $Z \sim \mathcal{P}2(m, h_1, \ldots, h_N)$ when it has density  
\begin{equation}
\label{Paretodensitymultivariate}
g_{m,h_1,\ldots, h_N}(z) \, = \, \frac{(m)_N \prod_{i=1}^N h_i}{(1+ \sum_{i=1}^N h_i z_i)^{m+N}} \, \mathbb{I}_{\mathbb{R}_+^N}(z),
\end{equation}
for parameters $h_1, \ldots h_N,m>0$, and $(m)_N= \frac{\Gamma(m+N)}{\Gamma(m)}$.
\end{definition} 

Such distributions form a $N+1$ parameter family with scale parameters $1/h_i, i=1, \ldots, N$, and shape parameter $m$.   As recorded with the following lemma containing a selection of properties that are known and readily verified, such multivariate Pareto distributions possess univariate Pareto marginals and conditionals, subvectors that are also distributed as multivariate Pareto with densities, and a joint survival function having a rather simple form. 

\begin{lemma}
\label{lemmapropertiesbivariatePareto}
Consider $Z \sim \mathcal{P}2(m, h_1, \ldots, h_N)$.   Let $i \neq j$.     Then, we have:  {\bf (i)}  $Z_i \sim f_{m,h_i}$, {\bf (ii)}  $Z_j|Z_i=z_i \sim f_{m+1, h_j/(1+h_iz_1)}$, {\bf (iii)}  $(h_1Z_i, \ldots, h_N Z_N)$ has density $g_{m,1, \ldots, 1}$,  {\bf (iv)} the joint survival function is given by $\mathbb{P}(\cap_{i=1}^{N} \{Z_i \geq z_i\}) \, = \, \{(1+ \sum_{i=1}^N  h_i z_i)\}^{-m}$ for $z_i \geq 0, i=1, \ldots, N$, {\bf (v)} the simple regressions are linear with {\bf (vi)} $\mathbb{E}(Z_j|Z_i) \, = \, \frac{1+h_iZ_i}{h_jm}\,$, and {\bf (vi)} the correlation between $Z_i$ and $Z_j$ is given by $\rho(Z_i, Z_j) = \frac{1}{m}$  for $m>2$.
\end{lemma}
 Furthermore,    we will require the following.

\begin{lemma}
\label{lemmasumparetos}
Let $Z \sim \mathcal{P}2(m, h_1, \ldots, h_N)$ and $W= \sum_{i=1}^N  h_i Z_i$.   Then $W$ is Beta type-II distributed (denoted $W \sim B2(N,m)$) with p.d.f.  $f_W(w) \, = \,  \frac{\Gamma(N+m)}{\Gamma(N) \Gamma(m)}  \frac{w^{N-1}}{(1+w)^{N+m}}$, for $w \in (0,\infty)$. 
\end{lemma}
{\bf Proof.}  The result is known but we provide a proof for completeness.  Given the multivariate Pareto representation $(h_1Z_1 \,\ldots, h_N Z_N) =^d   (\frac{E_1}{G}, \ldots, \frac{E_N}{G})$ with $E_1, \ldots E_N, G$ independently distributed, and with $E_i \sim Exp(1)$ and $G \sim \mathcal{G}(m,1)$, we see that $W$ is distributed as the ratio of two independent $\mathcal{G}(N,1)$ and $\mathcal{G}(m,1)$ variables, hence Beta type II with the given parameters.  \qed

\section{Predictive densities and regions}

Based on $X$ as in (\ref{modelsufficient}), we provide in this section Bayesian predictive densities and regions for $Y$ and $Z$ as defined in (\ref{predicted}) and (\ref{predictedmultivariate}).    We consider Gamma $\mathcal{G}(\alpha, \beta)$ prior densities  $\pi_{\alpha, \beta}(\theta) \propto \theta^{\alpha-1} \, e^{-\beta \theta} \, \mathbb{I}_{(0,\infty)}(\theta)$, including the usual non-informative case $\pi_0(\theta) \,=\, \frac{1}{\theta} \, \mathbb{I}_{(0,\infty)}(\theta)$ for $\alpha=\beta=0$.   
 
\subsection{Predictive densities}

We begin with the future next $N$ order statistics.   

\begin{theorem}
\label{theoremsimplecasemultivariate}
The Bayes predictive density of $Z$ in (\ref{predictedmultivariate}), based on $X \sim \mathcal{G}(m,\theta)$ and prior density $\pi_{\alpha, \beta}$ for $\theta$, is that of 
$\mathcal{P}2(m+\alpha, h_1, \ldots, h_N)$ distribution with $h_i= \frac{n-m-i+1}{x+\beta}$ for $i=1, \ldots, N$.
\end{theorem}
{\bf Proof.}  With $q_{\theta}(z|x) \,=\, \frac{ \theta^N \, (n-m)!}{(n-m-N)!} \, e^{- \theta \sum_{i=1}^N (n-m-i+1) z_i}$ and $\theta|x \sim \mathcal{G}(\alpha+m, x+\beta)$, we obtain from (\ref{predictivedensitygeneral}):
\begin{eqnarray*}
\nonumber  \hat{q}_{\pi_{\alpha,\beta}}(z|x)\, & = &\, \frac{(n-m)! \, (\beta + x)^{\alpha+m}} {(n-m-N)! \, \Gamma(\alpha+m) \, } \int_0^{\infty}  \theta^{N+\alpha+m-1}  \, e^{-\theta \big(\beta + x + \sum_{i=1}^N (n-m-i+1) z_i\big)}  \, d\theta  \, \\
\, & = & \,   \frac{(n-m)! \, (\beta + x)^{\alpha+m}} {(n-m-N)! \, \Gamma(\alpha+m) \, } \, \Gamma(N+\alpha+m) \, \big\{ \beta+ x + \sum_{i=1}^N (n-m-i+1) z_i \big\}^{-(N+\alpha+m)},
\end{eqnarray*}
which is indeed a $\mathcal{P}2(m+\alpha, h_1, \ldots, h_N)$ density.  \qed

Observe that the density has univariate Pareto marginals and multivariate Pareto distributed subvectors in accordance to Lemma \ref{lemmapropertiesbivariatePareto}.  

We now turn to the Bayesian predictive density for two future order statistics and demonstrate a non-probabilistic mixture of bivariate Pareto densities as given in Definition \ref{defmultivariatePareto} for $N=2$.

\begin{theorem}
\label{theoreompredictivedensity}
The Bayes predictive density of $Y$ in (\ref{predicted}), based on $X \sim \mathcal{G}(m,\theta)$ and prior density $\pi_{\alpha, \beta}$ for $\theta$, is given by
\begin{equation}
\label{expressionpredictivedensity}
\hat{q}_{\pi_{\alpha,\beta}}(y|x) \, = \,   \sum_{i=0}^{r-m-1} 
\sum_{j=0}^{s-r-1}  \omega_{i,j} \, g_{m+\alpha, \frac{a_i}{x+\beta}, \frac{b_j}{x+\beta}}  \big(y \big)\,,
\end{equation} 
with $g$ a bivariate Pareto density given in (\ref{Paretodensitymultivariate}) for $N=2$, 
$$a_i=n-r+i+1, b_j=n-s+j+1, \hbox{and } \omega_{i,j} \, = \, \frac{(n-m)!}{(n-s)!} \, \frac{(-1)^{i+j}}{i! \, j!} \, \,  \frac{1}{(r-m-i-1)! \, (s-r-j-1)!} \frac{1}{a_i b_j}.$$
\end{theorem}
{\bf Proof.}     From (\ref{lemmadensityofY}) and (\ref{predictivedensitygeneral}), with $\theta|x \sim \mathcal{G}(m+\alpha, x+\beta)$, we obtain
\begin{equation}
\nonumber
\hat{q}_{\pi_{\alpha,\beta}}(y|x) \, = \, \frac{(n-m)!}{(n-s)!} \frac{\,(x+\beta)^{m+\alpha}}{\Gamma(m+\alpha)} \, \int_0^{\infty}  \theta^{m+\alpha+1} \,  \frac{(1-e^{-\theta y_1})^{r-m-1}}{(r-m-1)!} \, \frac{(1-e^{-\theta y_2})^{s-r-1}}{(s-r-1)!}  e^{-\theta(L(y)+ x+ \beta)} \, d\theta,
\end{equation}
with $L(y)\,= \, (n-r+1) y_1 + (n-s+1) y_2$.   Binomial expansions and an interchange of sum and integral yield

\begin{eqnarray*}
\hat{q}_{\pi_{\alpha,\beta}}(y|x) \, & = &\, \frac{\,(x+\beta)^{m+\alpha}}{\Gamma(m+\alpha)} \, \sum_{i=0}^{r-m-1} 
\sum_{j=0}^{s-r-1} a_i b_j \, \omega_{i,j} \,  \int_0^{\infty}  \theta^{m+\alpha+1} \, e^{-\theta(L(y)+ i y_1 + j y_2 +  x+ \beta)} \, d\theta\, \\
& = &\,  (m+\alpha) \, (m+\alpha+1) \, (x+\beta)^{m+\alpha} \, \sum_{i=0}^{r-m-1} \sum_{j=0}^{s-r-1} a_i b_j \, \omega_{i,j} \, \big(L(y)+ i y_1 + j y_2 +  x+ \beta \big)^{-(m+\alpha+2)}, 
\end{eqnarray*}
which leads to (\ref{expressionpredictivedensity}).   \qed

We point out the following result for the next two order statistics, which follows immediately from either Theorem \ref{theoremsimplecasemultivariate} or Theorem \ref{theoreompredictivedensity}.
\begin{corollary}
\label{corollarysimplecase}
For the particular case where $r=m+1, s=m+2$,  the predictive density of $Y$ as in 
(\ref{predicted}) is that of a $\mathcal{P}2(m+\alpha, h_1, h_2)$ distribution
(\ref{expressionpredictivedensity}) with $h_1= \frac{n-m}{x+\beta}$ and $h_2= \frac{n-m-1}{x+\beta}$.
\end{corollary}

\begin{remark}
\label{remarkweights}
 Interestingly, the above weights $\omega_{i,j}$ arise through a series expansion of the Beta function.  Indeed, by the the binomial expansion of $(1-t)^{d-1}$ for $d \in \mathbb{N}_+$, we obtain for $c>0$
\begin{eqnarray*}
\int_0^1 t^{c-1} \, (1-t)^{d-1} \, dt \, & = & \, \frac{\Gamma(c) \Gamma(d)}{\Gamma(c+d)} \\
\, & \Longrightarrow   & \sum_{k=0}^{d-1} \gamma_{c,d,k}  =1 \,,
\end{eqnarray*} 
with $\gamma_{c,d,k} \, = \, \frac{\Gamma(c+d)}{\Gamma(c)} \,  \frac{(-1)^k}{k! \, (d-1-k)!} \,  \frac{1}{c+k}$.    Furthermore, observe that $\omega_{i,j}= \omega_{1,i}\, \omega_{2,j}$ with $\omega_{1,i}= \frac{(n-r)!}{(n-s)!} \,\gamma_{n-r+1, r-m, i}$ and $\omega_{2,j}= \frac{(n-s)!}{(n-r)!} \, \gamma_{n-s+1, s-r, j}$, demonstrating alternatively that $\sum_{i,j} \omega_{i,j}=1$.
\end{remark}

The  mixture representation of the predictive density in Theorem \ref{theoreompredictivedensity} coupled with the bivariate Pareto properties of Lemma \ref{lemmapropertiesbivariatePareto} facilitate the evaluation of the marginal and conditional distributions associated with (\ref{expressionpredictivedensity}).   As shown below, mixture representations of univariate Pareto distributions arise.

\begin{corollary}
\label{corollarymarginals}  
For $m < r \leq n-1$, the marginal density of $Y_1$ associated with the Bayes predictive density (\ref{expressionpredictivedensity}) is given by 
\begin{equation}
\label{densitymarginal}
\hat{q}_{\pi}(y_1|x)\, = \, \sum_{i=0}^{r-m-1} \gamma_{n-r+1, r-m, i} \, f_{m+\alpha, \frac{a_i}{x+\beta}}(y_1)\,,  
\end{equation}
where $\gamma_{c,d,k}$ is given in Remark \ref{remarkweights} and $f_{l,h}$ is a univariate Pareto density as given in (\ref{densityunivariatePareto}). In the particular case where $r=m+1$, we have that $\hat{q}_{\pi}(y_1;x) = f_{m+\alpha, \frac{n-m}{x+\beta}}(y_1)$ for $y_1>0$. 
\end{corollary}
{\bf Proof.}    We have from Theorem \ref{theoreompredictivedensity}
\begin{eqnarray*}
\hat{q}_{\pi}(y_1|x) \, & = & \, \int_0^{\infty} \hat{q}_{\pi}(y|x) \, dy_2 \\
\, & = & \,   \sum_{i=0}^{r-m-1} \gamma_{n-r+1, r-m, i} \sum_{j=0}^{s-r-1} \gamma_{n-s+1, s-r, j} \;  \int_{0}^{\infty} g_{m+\alpha, \frac{a_i}{x+\beta}, \frac{b_j}{x+\beta}}  \big(y \big) \, dy_2,
\end{eqnarray*}
which leads to (\ref{densitymarginal}) since the joint density $g_{m+\alpha, \frac{a_i}{x+\beta}, \frac{b_j}{x+\beta}}$ for $Y$ has marginal $f_{m+\alpha, \frac{a_i}{x+\beta}}$ for $Y_1$, for all $i,j$ (Lemma \ref{lemmapropertiesbivariatePareto}) and since $\sum_{j=0}^{s-r-1} \gamma_{n-s+1, s-r, j}=1$.  \qed
 
The result in itself is not new and was obtained with the univariate analysis carried out by  \cite{Dunsmore-1974}. A similar development establishes that (\ref{densitymarginal}) holds for the marginal predictive density of $X_{n:n}-X_{m:n}$.  
For the conditional distributions, we have the following.

\begin{corollary}
\label{corollaryconditionals}
For the Bayes predictive density (\ref{expressionpredictivedensity}): 
\begin{enumerate}
\item[\bf{ (a)}]   The conditional density of $Y_2$ given $Y_1=y_1$ is given by the mixture representation
\begin{equation}
\nonumber
\hat{q}_{\pi}(y_2|y_1;x) \, = \,  \sum_{i=0}^{r-m-1} \sum_{j=0}^{s-r-1} \alpha_i(y_1) \, \beta_j \, f_{m+\alpha+1, \frac{b_j}{x\, +\, \beta\,+\,a_i \,y_1}}
(y_2),  
\end{equation}
with $\alpha_i(y_1) \propto \gamma_{n-r+1, r-m, i} \, f_{m+\alpha, \frac{a_i}{x+\beta}}(y_1) $ such that $\sum_{i=0}^{r-m-1} \alpha_i(y_1)=1$, and $\beta_j= \gamma_{n-s+1, s-r, j}\,$;
\item[\bf{ (b)}]  The conditional density of $Y_1$ given $Y_2=y_2$ is given by the mixture representation 
\begin{equation}
\nonumber
\hat{q}_{\pi}(y_1|y_2;x) \, = \,  \sum_{i=0}^{r-m-1} \sum_{j=0}^{s-r-1} \xi_j(y_2) \, \nu_i \, f_{m+\alpha+1, \frac{a_i}{x\, +\, \beta\,+\,\beta_j \,y_2}}(y_1),  
\end{equation}
with $\xi_j(y_2) \propto \gamma_{n-s+1, s-r, j} \, f_{m+\alpha, \frac{b_j}{x+\beta}}(y_2) $  such that $\sum_{j=0}^{s-r-1} \xi_j(y_2)=1$, and $\nu_i= \gamma_{n-r+1,r-m, i}$.  
\end{enumerate}
\end{corollary}
{\bf Proof.}   Part {\bf (a)} follows directly with the properties of Lemma \ref{lemmapropertiesbivariatePareto} upon writing $  g_{m+\alpha, \frac{a_i}{x+\beta}, \frac{b_j}{x+\beta}}  \big(y \big) \, =  f_{m+\alpha, \frac{a_i}{x\, +\,  \beta }}
(y_1) \, f_{m+\alpha+1, \frac{b_j}{x\, +\, \beta\,+\,a_i y_1 }}
(y_2)$. The same approach leads to {\bf (b)}. \qed

\begin{remark}
\label{remarkmoments}
With the above marginal and conditional distributions expressible as mixtures of univariate Pareto distributions, corresponding moments are readily available.   For instance, one obtains
\begin{eqnarray*}
\mu_1=\mathbb{E}_{\pi}(Y_1|x) & = & \frac{x+\beta}{m+\alpha-1} \, \sum_{i=0}^{r-m-1} \gamma_{n-r+1,r-m,i} \, \frac{1}{a_i} \, \\
\hbox{ and } \mu_2(y_1) \, =  \, \mathbb{E}_{\pi}(Y_2|Y_1;x) \, & = & \, \sum_{i=0}^{r-m-1} \sum_{j=0}^{s-r-1} \alpha_i(y_1) \, \beta_j \, \frac{x+\beta + a_i y_1}{(m+\alpha) b_j}  \\ \, & = & \, y_1 \, \sum_{i=0}^{r-m-1} \sum_{j=0}^{s-r-1} \alpha_i(y_1) \, \beta_j \, \frac{a_i }{(m+\alpha) b_j} \, + \, \frac{x+\beta}{m+\alpha} \sum_{j=0}^{s-r-1} \,  \, \frac{\beta_j}{b_j}.
\end{eqnarray*}

Observe that the conditional expectation $\mathbb{E}_{\pi}(Y_2|y_1;x)$ is affine linear of the form $C+Dy_1$.
\end{remark}

\subsection{HPD Bayesian prediction regions}

A Bayesian prediction region with credibility $1- \lambda$ for $Z|\theta$ based on an observed value $x$ of $X|\theta$ and for a given prior density $\pi$ is such that 
\begin{equation}
\label{credibleregion}
\int_{R(x)} \hat{q}_{\pi}(z|x) \, dz \, = \, 1- \lambda,
\end{equation}
with $\hat{q}_{\pi}$ given in (\ref{predictivedensitygeneral}).  One such choice which minimizes volume is the ubiquitous HPD region which is of the form
\begin{equation}
\label{HPDform}
R_{HPD}(x) \, = \, \big\{z \in \mathbb{R}_+^N \,: \, \hat{q}_{\pi}(z|x) \, \geq k    \big\},
\end{equation}
where $k$ is chosen so that (\ref{credibleregion}) is satisfied.  Here is an explicit form of the HPD credible region for the future order statistics spacings $Z_1, \dots,  Z_N$.

\begin{theorem}
\label{theoremHPDmultivariate}
Based on $X_{i:n}, i=1, \ldots, m$, the first $m$ order statistics among $n$ from i.i.d. Exp$(\theta)$ data, setting $X$ as in (\ref{modelsufficient}) and $Z$  as in (\ref{predictedmultivariate}), the HPD region of credibility $1- \lambda$ for $Z$ associated with Gamma prior $\pi_{\alpha, \beta}$ for $\theta$ is given by
\begin{equation}
\label{HPDregionmultivariate}
R_{HPD}(x) \, = \, \big\{z \in \mathbb{R}_+^N \,: \, \sum_{i=1}^N (n-m-i+1) \, z_i   \, \leq c_0  (x+\beta)   \big\},
\end{equation}
with $c_0$ the quantile of order $1-\lambda$ of a $B2(N,m+\alpha)$ distribution.   
\end{theorem}
{\bf Proof.}   It follows from Theorem \ref{theoremsimplecasemultivariate} that $R_{HPD}$ defined generally in (\ref{HPDform}) is of form (\ref{HPDregionmultivariate}).
From Lemma \ref{lemmasumparetos}, the posterior predictive distribution of $W=\sum_{i=1}^N \frac{(n-m-i+1) \, z_i}{x+\beta}$   \, is $B2(N,m+\alpha)$ distributed and the result follows by setting $c_0$ such that $\mathbb{P}(W \leq c_0) = 1 -\lambda$. \qed

\begin{remark}   With simple integral or finite sum forms for the $B2(N,m+\alpha)$ c.d.f., the evaluation of the above quantiles is rather straightforward.  For instance in the bivariate case, we have $\mathbb{P}(W \leq c) = 1 - \frac{1 + c \, (m+\alpha+1)}{(1+c)^{m+\alpha+1}} $, so that $c_0$ is expressible as the solution in $c>0$ of
\begin{equation}
\nonumber
\label{c_0}
\frac{1 + c \, (m+\alpha+1)}{(1+c)^{m+\alpha+1}} \, = \, \lambda.
\end{equation}
\end{remark}

\begin{remark}
HPD credible regions are in general not invariant with respect to transformations, but they are whenever the transformation is affine linear.   Since the original order statistics $T=(X_{m+1:n}, \ldots, X_{m+N;n})^{\top}$ are related to the spacings $Z=(Z_1, \ldots, Z_N)^{\top}$ by the affine linear transformation  $T=b+AZ$, with $b=(X_{m:n}, \ldots, X_{m:n})^{\top}$ and $A$ a lower triangular matrix with non-zero elements $a_{i,j}= 1$ for $i \leq j$, the transformation of $R_{HPD}(x)$ given by (\ref{HPDregionmultivariate}) to the order statistics $X_{m:n}, \ldots, X_{m+N;n}$ is also HPD.

\end{remark}

 \subsection{Prediction regions:  an algorithm for the bivariate case}
\label{sectionalgorithm}

We present here a general bivariate case solution to obtain a prediction region with a given credibility based on Theorem \ref{theoreompredictivedensity}'s predictive density for $Y$ as in (\ref{predicted}).  There are several options available, but we opt for a rather explicit strategy to construct a Bayesian credible region.
It can be applied for any choice of $(r,s,m,n)$, prior $\pi_{\alpha, \beta}$ such that $m+ \alpha>1$, and the targeted credibility $1-\lambda$.    We refer to notation used throughout this paper, namely $a_i$ and $b_j$ as in Theorem \ref{theoreompredictivedensity}, $\gamma_{c,d,k}$ as in Remark \ref{remarkweights}, $f_{l,h}$ and $\bar{F}_{l,h}$ as the univariate Pareto density and survival functions in (\ref{densityunivariatePareto}), $\beta_j$ and $\alpha_i(y_1)$ as in Corollary \ref{corollaryconditionals}, and $\mu_1$ and $\mu_2(y_1)$ as in Remark \ref{remarkmoments}.

\begin{enumerate}
\item[]Step 1.  Select a prediction region $A$ for $Y_1$ of credibility $\sqrt{1-\lambda}$
based on the marginal density in Corollary \ref{corollarymarginals}.    Such a choice would desirably be an interval with relatively high levels of the predictive density.  A suitable choice is:
\begin{equation}
\nonumber
 A\, = \, \big[\mu_1- \Delta_1, \mu_1+ \Delta_1\big] \, \cap [0,\infty) \,,
\end{equation}
with $\Delta_1 >0$ uniquely chosen such that
\begin{eqnarray*}
& & \int_{\mu_1-\Delta_1}^{\mu_1+\Delta_1}  \hat{q}_{\pi_{\alpha, \beta}}(y_1|x) \, dy_1 \, =  \, \sqrt{1 - \lambda} \\
\Longleftrightarrow
& & \sum_{i=0}^{r-m-1} \gamma_{n-r+1,r-m,i} \, \Big(\bar{F}_{m+\alpha, \frac{a_i}{x+\beta}}(\mu_1-\Delta_1)_+ \, - \,  \bar{F}_{m+\alpha, \frac{a_i}{x+\beta}}(\mu_1+\Delta_1)    \Big) \, = \, \sqrt{1 - \lambda},
\end{eqnarray*}
with $z_+=\max\{0,z\}$.   Since this last expression is strictly increasing in $\Delta_1$, it is rather straightforward to approach $\Delta_1$ numerically. 

\item[]Step 2.  For each $y_1 \in A$, select a prediction region $B(y_1)$ of conditional credibility $\sqrt{1-\lambda}$ for $y_2$ based on the conditional density of $Y_2|Y_1=y_1$ given in part (a) of Corollary \ref{corollaryconditionals}. The challenge here is similar to the one in Step 1 but to be repeated for all $y_1$.   Analogously to Step 1, $B(y_1)$ can reasonably be constructed pivoting around the mean as 
\begin{equation}
\nonumber B(y_1) \, = \, \big[\mu_2(y_1)- \Delta_2(y_1), \mu_2(y_1)+ \Delta_2(y_1)\big] \, \cap [0,\infty),
\end{equation}
with $\Delta_2(y_1) > 0$  chosen such that

\begin{eqnarray*}
& & \! \sum_{i=0}^{r-m-1} \sum_{j=0}^{s-r-1} \alpha_i(y_1) \, \beta_j \, \int_{\mu_2(y_1)-\Delta_2(y_1)}^{\mu_2(y_1)+\Delta_2(y_1)} \, f_{m+\alpha+1, \frac{b_j}{x+\beta+a_i y_1}}(y_2) \, dy_2 \, =  \, \sqrt{1 - \lambda} \\
\Longleftrightarrow    \!
&  & \sum_{i=0}^{r-m-1} \sum_{j=0}^{s-r-1} \alpha_i(y_1) \, \beta_j \, \, \Big(\bar{F}_{m+\alpha+1, \frac{b_j}{x+\beta+a_i y_1}}(\mu_2(y_1)-\Delta_2(y_1))_+ \, - \,  \bar{F}_{m+\alpha+1, \frac{b_j}{x+\beta+a_i y_1}}(\mu_2(y_1)+\Delta_2(y_1))    \Big) \,  \, \\ & & = \sqrt{1 - \lambda}.
\end{eqnarray*}

\end{enumerate}

The resulting prediction region $R=\{ (y_1, y_2):   y_1 \in A, y_2 \in B(y_1) \}$ has credibility $1-\lambda$ indeed since 
\begin{eqnarray*}
\mathbb{P}_{\pi}\big(Y_1 \in A, Y_2 \in B(y_1)|x  \big) & = & \int_{A} \hat{q}_{\pi}(y_1|x) \big\{ \int_{B(y_1)} \hat{q}_{\pi}(y_2|y_1;x) \, 
  dy_2 \big\} \,  dy_1 \\
\, & = & \int_{A} \hat{q}_{\pi}(y_1|x) \, (\sqrt{1-\lambda}) \, dy_1 \\ 
& = &  1 - \lambda.
\end{eqnarray*}

\begin{remark}
Other regions can also be selected.  For instance, one-sided choices with $A$ of the form $[0,\bar{\Delta}_1]$ or $[\underline{\Delta}_1, \infty)$, and with $B$ of the form $[0,\bar{\Delta}_2(y_1)]$ or $[\underline{\Delta}_2(y_1), \infty)$.  Another alternative would be to aim for different credibilities $1-\lambda_1$ and $1-\lambda_2$ in Steps 1 and 2, respectively, such that  $(1-\lambda_1)(1-\lambda_2)=1-\lambda$.
\end{remark}

\subsection{Example}

The following dataset from \cite{Murthy-et-al} shows $n=30$ ordered failure times for repairable items:
\begin{table*}[h!]
	\centering
	\begin{tabular}{cccccccccc}
		\hline
		0.11& 0.30& 0.40& 0.45& 0.59& 0.63& 0.70& 0.71& 0.74& 0.77\\
		0.94& 1.06& 1.17& 1.23& 1.23& 1.24& 1.43& 1.46& 1.49& 1.74\\
		1.82& 1.86& 1.97& 2.23& 2.37& 2.46& 2.63& 3.46& 4.36& 4.73\\
		\hline
	\end{tabular}	
\end{table*}

For the purpose of illustration, suppose that only the $m=20$ first values are observed and that we wish to predict either:
\begin{enumerate}
\item[ {\bf (i)}]  the next two order statistics spacings $(Z_1, Z_2) \,=\, \big(X_{21:30}-X_{20:30}, X_{22:30}-X_{21:30}\big)$,
\end{enumerate}
or
\begin{enumerate}
\item[ {\bf (ii)}]  the next and last order statistics  $\big(X_{21:30}, X_{30:30}\big)$.
\end{enumerate}
In both cases, we use the non-informative prior density $\pi_(\theta)\,=\, \frac{1}{\theta} \, \mathbb{I}_{(0,\infty)}(\theta)$, and consider credibility $1-\lambda=0.95$.  As expanded upon in the next section, the frequentist coverage of such a prediction region matches the credibility for all $\theta>0$.  The data yields $x\,=\, \sum_{i=1}^{20} x_{i:30} \, + \, 10 \, x_{20:30}\, = \, 35.79$.

For case {\bf (i)}, Theorem \ref{theoremHPDmultivariate} applies and the HPD credible region for $(Z_1, Z_2)$ is given by 
\begin{align*}
R_{\textrm{HPD}}(35.79) \, = \, \bigl\{(z_1,z_2) \in \mathbb{R}_+^2 \,:\, 0.2794 \, z_1\, + \, 0.2515 z_2 \, \leq 0.2606\big\},
\end{align*}
$c_0\,=\,0.2606$ being the quantile of order $0.95$ of a $B2(2,20)$ distribution.

\begin{figure}[ht!]
	\centering	
\subfigure[$1-\lambda=0.95$]{
	\includegraphics[scale=0.4]{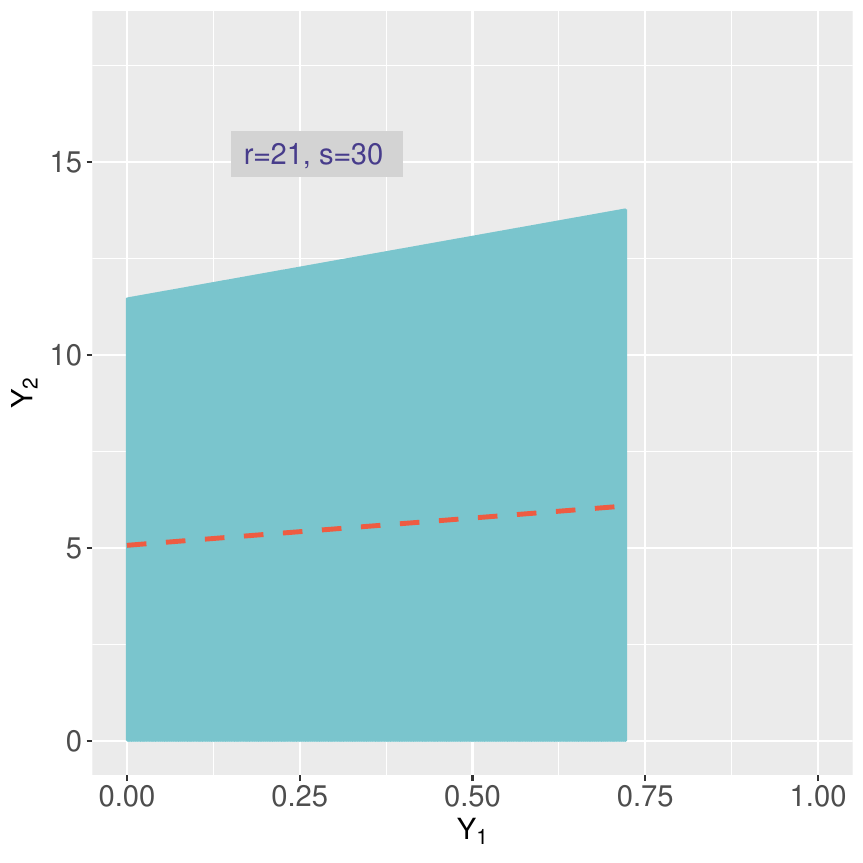}}
	\qquad
\subfigure[$1-\lambda=0.80$]{
	\includegraphics[scale=0.4]{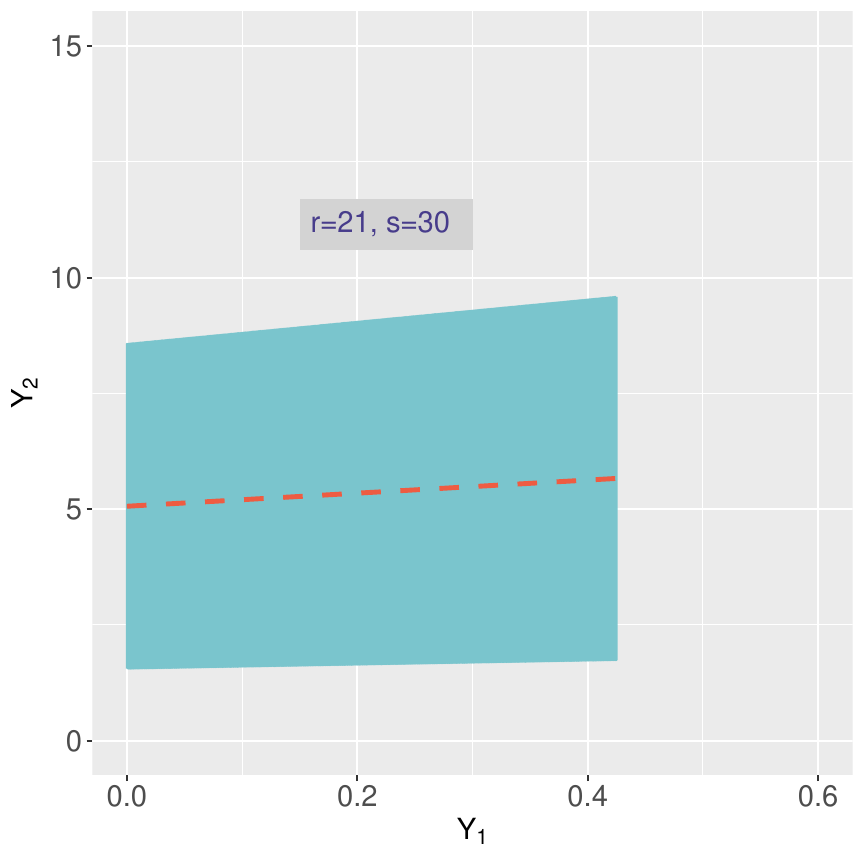}}	
	\caption{Bayesian prediction regions and $\mathbb{E}(Y_2|y_1,x=35.79)$ (dashed) } \label{fig-region}
\end{figure}

\begin{figure}[ht!]
	\centering	
	\subfigure[$1-\lambda=0.95$]{
		\includegraphics[scale=0.4]{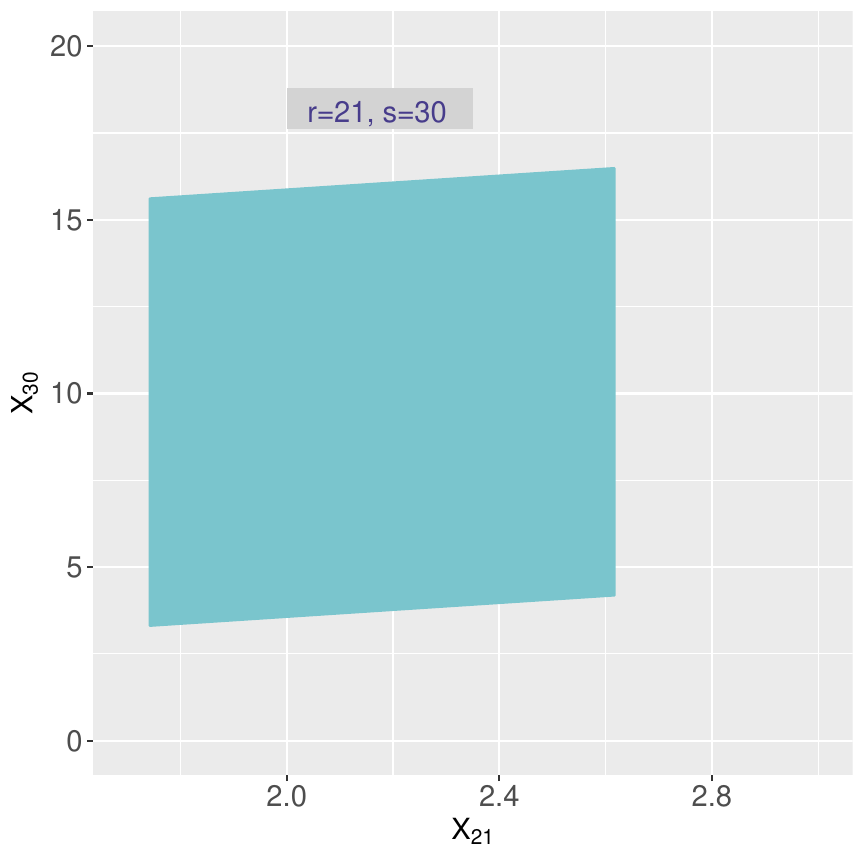}}
	\qquad
	\subfigure[$1-\lambda=0.80$]{
		\includegraphics[scale=0.4]{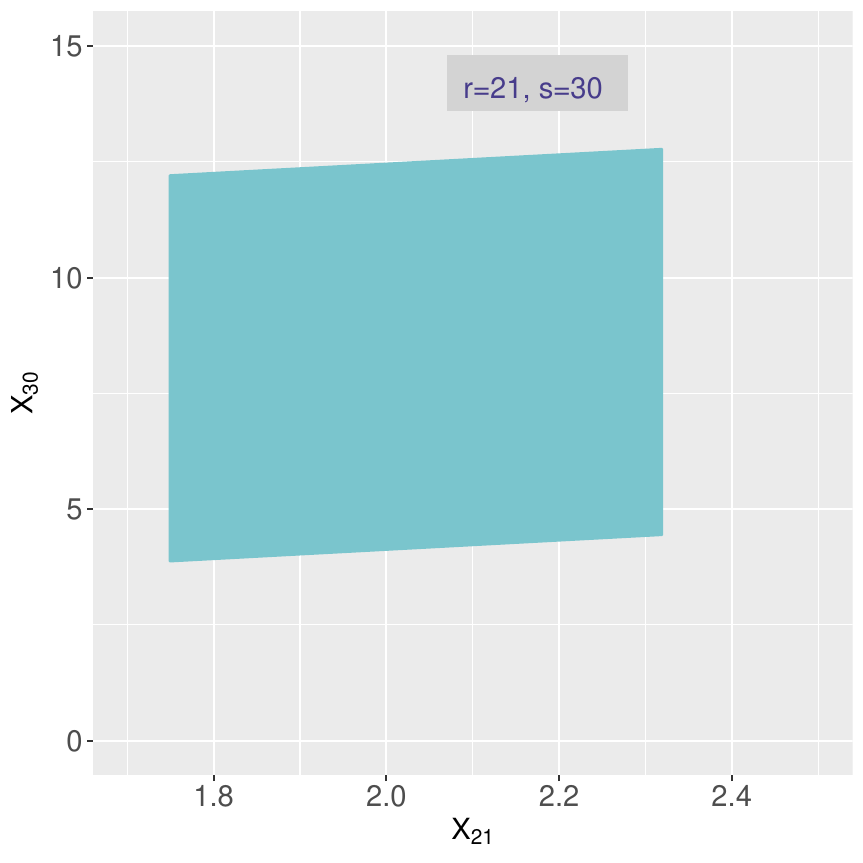}}	
	\caption{
			Bayesian prediction regions for $(X_{21:30},X_{30:30})$ in case (ii)} \label{fig-region2}
\end{figure}

For case {\bf (ii)}, we illustrate the use of Section \ref{sectionalgorithm}'s algorithm which is prescribed for $(Y_1, Y_2)$ with $Y_1 \, = \, X_{21:30}-X_{20:30}$ and $Y_2 \,=  X_{30:30}-X_{21:30}$.  Figure \ref{fig-region} (a) presents the resulting Bayesian prediction region along with the regression function or conditional expectation $\mathbb{E}(Y_2|y_1,x=35.79)$ as a function of $y_1$.  The first step yields $A=[0, 0.722]$, while the $B(y_1)$ intervals are shown in the figure. For example, $B(0.50)=\, [0, 13.059]$. 
Observe that a left neighbourhood of $0$ is included in $B(y_1)$ for all $y_1 \in A$ as the second step credibility, equal to $(0.95)^{1/2} \approx 0.9747$, is quite large. In contrast, the Bayesian prediction region with credibility $0.80$ is displayed in Figure \ref{fig-region} (b) resulting in $A=[0,0.4258]$ and intervals $B(y_1)$ centered at $\mathbb{E}(Y_2|y_1,x=35.79)$ which exclude values close to $0$ for all $y_1 \in A$.  Finally, the corresponding prediction regions for $\big(X_{21:30}, X_{30:30}\big) \,$ and credibilities $0.95$ and $0.80$ are obtained as $\{(y_1+x_{20:30}, y_1+y_2+x_{20:30}): (y_1, y_2) \in R \}$ and displayed in Figure \ref{fig-region2}. 

\section{Frequentist coverage and credibility}

Through Theorem \ref{theoreompredictivedensity}'s Bayesian predictive densities, the prediction regions for future order statistics described in the previous section are constructed in order to attain exact Bayesian credibility.  This includes the non-informative prior prediction densities $\hat{q}_{\pi_0}(\cdot;X)$ given by Theorems \ref{theoremsimplecasemultivariate} and  \ref{theoreompredictivedensity} for $\alpha=\beta=0$. Moreover, prediction regions based on $\hat{q}_{\pi_0}(\cdot;X)$  also lead to exact frequentist coverage as expanded on below.   In fact, we cast the result in a more general scale parameter family model setting with scale parameter densities 
\begin{equation}
\label{modelscale}
p_{\sigma}(x) = \frac{1}{\sigma} \, p_1(\frac{x}{\sigma}) \hbox{ and } q_{\sigma}(y') = \frac{1}{\sigma^{d_2}} \, q_1(\frac{y_1'}{\sigma}, \ldots, \frac{y_{d_2}'}{\sigma}).
\end{equation}

We will make use of the following intermediate result, a univariate version of which was first given by  \cite{LMoudden-et-al} (i.e., $d_2=1$).

\begin{lemma}
\label{lemmaratio}
Under model (\ref{modelscale}) and prior density $\pi_0(\sigma) \, = \, \frac{1}{\sigma} \, \mathbb{I}_{(0,\infty)}(\sigma)$, the Bayesian predictive density $\hat{q}_{\pi_0}$ is given by 
\begin{equation}
\label{predictivedensityratio}  \hat{q}_{\pi_0}(y'|x) \, = \, \frac{1}{x^{d_2}}  \, h(\frac{y'}{x}), 
\end{equation}
where $h$ is the frequentist density of $R=\frac{Y'}{X}=(R_1, \ldots, R_{d_2})$, which is free of $\sigma$, and given by  
\begin{equation}
\label{ratiodensity}
h(r) \, = \,  \int_0^{\infty} u^{d/2} \, q_1(r_1u, \ldots, r_{d_2} u) \, p(u) \, du.
\end{equation}
Furthermore, it is the case that 
\begin{equation}
\label{equalitydistribution}
\big(\frac{Y_1'}{x}, \ldots, \frac{Y'_{d_2}}{x}\big) \big| x  \, =^d \,  \big(\frac{Y_1'}{X}, \ldots, \frac{Y'_{d_2}}{X}\big) \big| \sigma   \, \hbox{ for all }   x, \sigma,
\end{equation}
i.e., the posterior predictive and frequentist distributions of $R$ match and are furthermore independent of the observed value of $x$ and the parameter value of $\sigma$. 
\end{lemma}

{\bf Proof.}   Identity (\ref{equalitydistribution}) follows from the first part of the lemma.
A direct evaluation yields (\ref{ratiodensity}) as the density of $R|\sigma$.  Finally, the posterior density of $\sigma$ is given by $\pi(\sigma|x) \, = \, \frac{x}{\sigma^2} \, p_1(\frac{x}{\sigma})$, which leads to the predictive density
\begin{eqnarray*}
\hat{q}_{\pi_0}(y'|x) \, &  = & \, \int_0^{\infty}  \frac{1}{\sigma^{d_2}} \, q_1(\frac{y_1'}{\sigma}, \ldots, \frac{y_{d_2}'}{\sigma}) \, \frac{x}{\sigma^2} \, p_1(\frac{x}{\sigma}) \, d\sigma \\
\, & = & \frac{1}{x^2}  \int_0^{\infty} u^{d_2}  q_1\big(\frac{y_1' u}{x}, \ldots, \frac{y_{d_2}' u }{x}   \big) \, p(u) \, du  \\
\, & = &  \frac{1}{x^{d_2}}  \, h(\frac{y'}{x})\,.  \;\;\; \; \; \; \;\;\;\;\;\; \; \; \; \;\;\; \;\;\; \; \; \; \;\;\;  \qed
\end{eqnarray*}

\begin{remark}
Identity (\ref{equalitydistribution}) is quite general and clarifies why Bayesian analysis with respect to the prior density $\pi_0$ matches pivotal based analysis that stems from the right-hand side.   In the univariate case and for the prediction of a single future order statistic, as remarked upon by Dunsmore (1974), the identity explains leads to his Bayesian solutions matching the pivotal based solution of Lawless (1971).
\end{remark}

\begin{theorem}
\label{coverage}
Consider model (\ref{modelscale}) and a Bayesian prediction region $R(X)$ with credibility $1- \lambda$ associated with the prior density $\pi_0(\theta)=\frac{1}{\theta} \, \mathbb{I}_{(0,\infty)}(\theta)$.    Then,  $R(X)$ has exact frequentist coverage probability, i.e., $\mathbb{P}(R(X) \ni Y'|\theta) \, = \,1- \lambda$ for all $\theta>0$.
\end{theorem}
{\bf Proof.}  Since $R(X)$ has credibility $1- \lambda$, we have 
\begin{equation}
\nonumber   \mathbb{P}\big(Y' \in R(x))|x\big) \, = \, \mathbb{P}\big( \frac{Y'}{x} \in R^*(x) \big|x\big) \, = \, 1- \lambda,
\end{equation}
where $R^*(x) \, = \, \{r \in \mathbb{R}_{+}^{d_2}: \frac{r}{x} \in R(x)\}$.  
Now, since the predictive distribution of $\frac{Y'}{x}$ is free of $x$ with p.d.f. $h$ (Lemma \ref{lemmaratio}), it follows that $R^*(x)$ is free of $x$ such that 
$\int_{R^*(x)} h(z) \, dz \, = \, 1- \lambda $.  
On the other hand, the frequentist coverage of $R(X)$ is given by
\begin{equation}
\nonumber
\mathbb{P}(Y' \in R(X)|\,\theta\,) \, = \, \mathbb{P}\big(\frac{Y'}{X} \in R^*(X)\big| \theta \big) \, = \, \int_{R^*(x)} h(r) \, dr \, = \, 1- \lambda,
\end{equation}
since $h$ is also the density of $\frac{Y'}{X}|\theta$ (Lemma \ref{lemmaratio}).  \qed

To conclude, the above (with $Y'=Y$ or $Y'=Z$) applies to our set-ups as follows.

\begin{corollary}
\label{corollarycoverage}
Based on $X$ as in (\ref{modelsufficient}) and the non-informative prior density 
$\pi_0(\theta)=\frac{1}{\theta} \, \mathbb{I}_{(0,\infty)}(\theta)$, Bayesian predictive regions for $Z$ or $Y$ with credibility $1-\lambda$, based on the predictive densities given in Theorems \ref{theoremsimplecasemultivariate} and \ref{theoreompredictivedensity}, have matching frequentist coverage probability $1-\lambda$ for all $\theta>0$.
\end{corollary}

\section{Best invariance and minimaxity}

As seen in Section 3, Bayesian predictive densities such as those given in Theorem \ref{theoreompredictivedensity} facilitate the construction of a prediction region for future values of $Y$ or $Z$ with a given credibility.  Non-Bayesian predictive densities are also available, such as plug-in densities and those obtained by likelihood or pivotal-based methods.   Kullback-Leibler (KL) divergence loss and  accompanying risk can be used to evaluate the frequentist performance of density estimators $\hat{q}(\cdot;X)$.  For our problem, these are given by
$$  L_{KL}(\theta, \hat{q}(\cdot;x)) \, = \, \int  q_{\theta}(t) \, \log \big(\frac{q_{\theta}(t)}{\hat{q}(t;x)}  \big) \, dt\,$$
(assuming Lebesgue densities), and
$$ R_{KL}(\theta, \hat{q}) \, = \,  \mathbb{E}_{\theta} \big\{ L_{KL}(\theta, \hat{q}(\cdot;X))  \big\}.$$
 It is of interest to assess the efficiency of predictive densities for KL risk, and decision-theoretic properties of invariance and minimaxity applicable to our contexts are reviewed in this section.    An early contribution to the determination of a best invariant density is due to  \cite{Murray-1977}, while a more exhaustive treatment of best invariant densities, as well as minimaxity in predictive density estimation, appears in \cite{Liang-Barron}.

A density $\hat{q}_m$ is minimax whenever $\hat{q}_m$ minimizes among all densities the supremum frequentist risk, i.e., in our cases when
\begin{equation}
\nonumber
\sup_{\theta > 0}  R_{KL}(\theta, \hat{q}_m) \, = \, \inf_{\hat{q}} \sup_{\theta >0}  R_{KL}(\theta, \hat{q})\,.
\end{equation}

For our problems, and more generally model (\ref{modelscale}) with $\sigma=1/\theta$, Kullback-Leibler divergence loss, a predictive density $\hat{q}$ is invariant under  changes of scale whenever it satisfies the scale parameter family property
\begin{equation}
\label{propertyinvariance}  
\hat{q}(y';c x) \, = \, \frac{1}{c^{d_2}}  \, \hat{q}(\frac
{y'}{c};x)\,, y' \in \mathbb{R}_+^{d_2} \,, 
\end{equation}
for all $c, x >0$.   The class of invariant densities here includes $\hat{q}_{\pi_0}$; as can be verified directly from the expressions given in Theorems \ref{theoremsimplecasemultivariate} and \ref{theoreompredictivedensity} with $\alpha=\beta=0$, or even by (\ref{predictivedensityratio}); as well as plug-in densities $q_{\hat{\theta}}$ with $\hat{\sigma}(x)=kx$, i.e., a scale invariant point estimator $\hat{\sigma}$ of $\theta$ satisfying  $\hat{\sigma}(c x) \, = \, c \, \hat{\theta}(x)$ for $c, x >0$; such as the maximum likelihood choice with $c=1/m$.   

The present invariance structure implies that an invariant density has constant risk as a function of $\theta$ as long as it is finite, from which it follows that there exists an optimal choice among invariant densities.  The risk constancy follows for vastly more general settings (e.g.,  Berger, 1985), but can be also derived directly as follows:
\begin{eqnarray*}
R_{KL}(c \sigma, \hat{q}) \, & = &  \int_{(0,\infty)} p_{c \sigma}(x) \, \int_{(0,\infty)^{d_2}}  q_{c \sigma}(y') \, \log\big(\frac{q_{c \sigma}(y')}{\hat{q}(y;x)}  \big) \, dy' \, dx \\
\, & = &  \int_{(0,\infty)} \frac{1}{c} \, p_{\sigma}(\frac{x}{c}) \, \int_{(0,\infty)^{d_2}}  q_{\sigma}(\frac{y'}{c}) \frac{1}{c^{d_2}} \, \log\big( \, \frac{q_{\sigma}(\frac{y'}{c})}{\hat{q}(\frac{y'}{c} ;\frac{x}{c} )}  \big) \, dy' \, dx \\
\, & = &  \int_{(0,\infty)} p_{\sigma}(t) \, \int_{(0,\infty)^{d_2}}  q_{\sigma}(u) \, \log\big(\frac{q_{\sigma}(u)}{\hat{q}(u;t)}  \big) \, du \, dt  \\   \, & = & \, R_{KL}(\sigma, \hat{q}),
\end{eqnarray*}
for $c, \sigma>0$, with the change of variables $(t,u) \, = \, \frac{1}{c} (x,y')$, and by making use of (\ref{propertyinvariance}).

Under general conditions for problems that are invariant, which are met here, a best invariant procedure exists and coincides with the generalized Bayes estimator associated with a (right) invariant prior density (e.g., Berger, 1985, section 6.6.2).  In out set-up, such a prior density is the non-informative $\pi_0$, and it leads to the best invariant property of the density $\hat{q}_{\pi_0}$.   Furthermore, $\hat{q}_{\pi_0}$ is minimax (see  \citet{Liang-Barron}, Theorem 2 and Proposition 3).  We conclude this section by summarizing the above as applicable (with $Y'=Y$ or $Y'=Z$) to the problems at hand.

\begin{theorem}
\label{theoreminvariance}
Consider $X_{i:n}, i=1, \ldots, m$, the first $m$ order statistics among $n$ from i.i.d. Exp$(\theta)$ data, and $X$, $Y$, and $Z$ as in (\ref{modelsufficient}),  (\ref{predicted}) and (\ref{predictedmultivariate}).  Then the best invariant densities for estimating the densities of $Y$ and $Z$, respectively, under KL divergence loss are $\hat{q}_{\pi_0}$, as given in Theorems \ref{theoremsimplecasemultivariate} and \ref{theoreompredictivedensity} for $\alpha=\beta=0$.   Furthermore, its KL risk is constant as a function of $\theta$ and $\hat{q}_{\pi_0}$ is minimax. 
\end{theorem}

\section{Concluding remarks}

We have illustrated the natural usage of Bayesian credible regions for the prediction of future order statistics under a type-2 censoring scheme with exponentially distributed data.   We have emphasized multivariate aspects and provided explicit expressions for Bayesian predictive densities and resulting prediction credible regions.  We have also addressed optimality properties achieved with the non-informative prior density choice, such as the matching of Bayesian credibility and frequentist probabilty coverage, as well as the best invariant and minimax properties of the corresponding Bayesian predictive density.

It would be interesting to extend the analysis more broadly to other types of probability models and censoring schemes.  We have provided a Bayesian HPD credible region (Theorem \ref{theoremHPDmultivariate}) for the next $N$ future order statistics.  However, for the general bivariate case, such a solution is lacking.   We do not know for instance if the Bayesian predictive density is unimodal which would facilitate the determination of such a region.

\section*{Acknowledgements}

The authors are thankful to the Universit\'e de Sherbrooke for its financial support
awarded as part of its visiting researcher's program.  \'Eric Marchand's research is supported in part by the NSERC of Canada.


\begin{thebibliography}{100}

\bibitem[Arnold (2014)]{Arnold-2014}  Arnold, B.C.   (2014).  Univariate and multivariate Pareto models.   {\it Journal of Statistical Distributions and Applications}, {\bf 1}, article 11.

\bibitem[Bagheri et al. (2022)]{Bagheri-et-al}  Bagheri, S.F., Asgharzadeh, A., Fern\'andez, A.J., \& P\'erez-Gonz\`alez, C.P. (2022).   Joint prediction for future failure times under type-II censoring.   {\it IEEE Transactions on Reliability}, {\bf 71}, 100--110.

\bibitem[Berger (1985)]{Berger-1985}  Berger, J.  (1985).   {\it  Statistical Decision Theory and Related Topics}.  Springer Texts in Statistics. Springer-Verlag, New York.  Second edition.   

\bibitem[Dunsmore (1974)]{Dunsmore-1974}  Dunsmore, I.R. (1974).   The Bayesian prediction problem in life testing models.   {\it Technometrics}, {\bf 16}, 455--4600.

\bibitem[Fourdrinier et al. (2019)]{Fourdrinier-et-al}   Fourdrinier, D.,  Marchand, É. \& Strawderman, W.E. (2019).  On efficient prediction and predictive density estimation for spherically symmetric models. {\it Journal of Multivariate Analysis}, {\bf 173}, 18--25.


\bibitem[Kotz et al. (2000)]{Kotz-et-al}  Kotz, S., Balakrishnan, N.  \& Johnson, N.L. (2000).   {\it Multivariate continuous distributions:  volume 1},  Wiley.

\bibitem[Lawless (1971)]{Lawless-1971}  Lawless,  J.F. (1971).   A prediction problem concerning samples from the exponential distributions, with application in life testing.   {\it Technometrics}, {\bf 13}, 725--730.

\bibitem[Lehmann \& Casella (1998)]{Lehmann-Casella} Lehmann, E.L. \& Casella, G.  (1998).  {\it Theory of Point Estimation}. Second edition,  Springer, New York.

\bibitem[Liang \& Barron (2004)]{Liang-Barron}
Liang,  F. \& Barron, A.  (2004). Exact minimax strategies for predictive density estimation, data compression, and model selection. 
{\it IEEE Transactions on Information Theory}, {\bf 62}, 2708--2726.


\bibitem[L'Moudden et al. (2017)]{LMoudden-et-al}  L'Moudden, A., Marchand, \'{E}., Kortbi, O. \& Strawderman, W.E. (2017).  On predictive density estimation for Gamma models with parametric constraints.  {\it Journal of Statistical Planning and Inference}, {\bf 185}, 56--68.

\bibitem[Murray (1977)]{Murray-1977}  Murray,  G.D. (1977).  A note on the estimation of probability density functions.  {\it Biometrika}, {\bf 64}, 150--152.

\bibitem[Murthy et al. (2004)]{Murthy-et-al}  Murthy, D.P., Xie, M. \& Jiang, R.  (2004).   {\it Weibull models}.   John Wiley \& Sons.   

\end{thebibliography}
\end{document}